\documentclass[conference]{IEEEtran}
\IEEEoverridecommandlockouts
\usepackage{cite}
\usepackage{amsmath,amssymb,amsfonts}
\usepackage{algorithmic}
\usepackage{graphicx}
\usepackage{textcomp}
\usepackage{xcolor}
\def\BibTeX{{\rm B\kern-.05em{\sc i\kern-.025em b}\kern-.08em
    T\kern-.1667em\lower.7ex\hbox{E}\kern-.125emX}}

\DeclareMathOperator*{\SubjectTo}{Subject\phantom{a}to:}
\DeclareMathOperator*{\Minimize}{Minimize:}

\newcommand{\norm}[1]{\left\lVert#1\right\rVert}
\newcommand{\xddots}{%
  \raise 4pt \hbox {.}
  \mkern 6mu
  \raise 1pt \hbox {.}
  \mkern 6mu
  \raise -2pt \hbox {.}
}

\begin{document}

\title{Battery Degradation Heuristics for Predictive Energy Management in Shipboard Power Systems \\
}

\author{\IEEEauthorblockN{1\textsuperscript{st} Satish Vedula}
\IEEEauthorblockA{\textit{Dept. of Electrical and Computer Engineering} \\
\textit{FAMU-FSU College of Engineering}\\
Tallahassee, USA \\
svedula@fsu.edu}
\vspace{2mm}
\IEEEauthorblockN{2\textsuperscript{nd} Ayobami Olajube}
\IEEEauthorblockA{\textit{Dept. of Electrical and Computer Engineering} \\
\textit{FAMU-FSU College of Engineering}\\
Tallahassee, USA \\
aolajube@fsu.edu}
\and
\IEEEauthorblockN{3\textsuperscript{rd} Koto Omiloli}
\IEEEauthorblockA{\textit{Dept. of Electrical and Computer Engineering} \\
\textit{FAMU-FSU College of Engineering}\\
Tallahassee, USA \\
kao23a@fsu.edu}
\vspace{2mm}
\centering
\IEEEauthorblockN{4\textsuperscript{th} Olugbenga Moses Anubi}
\IEEEauthorblockA{\textit{Dept. of Electrical and Computer Engineering} \\
\textit{FAMU-FSU College of Engineering}\\
Tallahassee, USA \\
oanubi@fsu.edu}
}

\maketitle

\begin{abstract}
The presence of Pulse Power Loads (PPLs) in the Notional Shipboard Power System (SPS) presents a challenge in the form of meeting their high ramp rate requirements. Considering the ramp rate limitations on the generators, this might hinder the power flow in the grid. Failure to meet the ramp rate requirements might cause instability. Aggregating generators with energy storage elements usually addresses the ramp requirements while ensuring the power demand is achieved. This paper proposes an energy management strategy that adaptively splits the power demand between the generators and the batteries while simultaneously considering the battery degradation and the generator's efficient operation. Since it is challenging to incorporate the battery degradation model directly into the optimization problem due to its complex structure and the degradation time scale which is not practical for real-time implementation, two reasonable heuristics in terms of minimizing the absolute battery power and minimizing the battery state of charge are proposed and compared to manage the battery degradation. A model predictive energy management strategy is then developed to coordinate the power split considering the generator efficiency and minimizing the battery degradation based on the two heuristic approaches. The designed strategy is tested via a simulation of a lumped notional shipboard power system. The results show the impact of the battery degradation heuristics for energy management strategy in mitigating battery degradation and its health management.
\end{abstract}


\section{Introduction}
In Shipboard Power Systems (SPSs), power generation modules (PGMs) provide the power required by various loads through the direct current (DC) distribution system. Modern SPSs are equipped with advanced power loads such as rail guns, electromagnetic radars, and heavy nonlinear and pulsed power loads (PPLs), these advanced loads are high ramp rate loads \cite{Derry_1}\cite{Derry_2}. SPSs consist of two key power-supplying sources namely: Power Generation Modules and Power Conversion Modules (PCMs) and power consuming loads: Power Load Modules (PLM). PGMs are ramp rate limited which hinders in meeting the high ramp rate power required by the loads. Failure to ramp up in time to meet the load requirement might lead to system imbalances and instability. Adding additional generators is not a feasible solution considering the restrictions on the ship hull dimensions. Thus, the PCMs, which can support high ramp rates provide a solution in addressing the load requirement of high ramp rate loads \cite{Kim_2015}\cite{ESRDC_1}. The PCMs can operate bi-directionally, they can store the energy and provide immediate support when high ramp rate loads are requesting power. 

\begin{figure}[t!]
      \centering
      \includegraphics[width=0.48\textwidth]{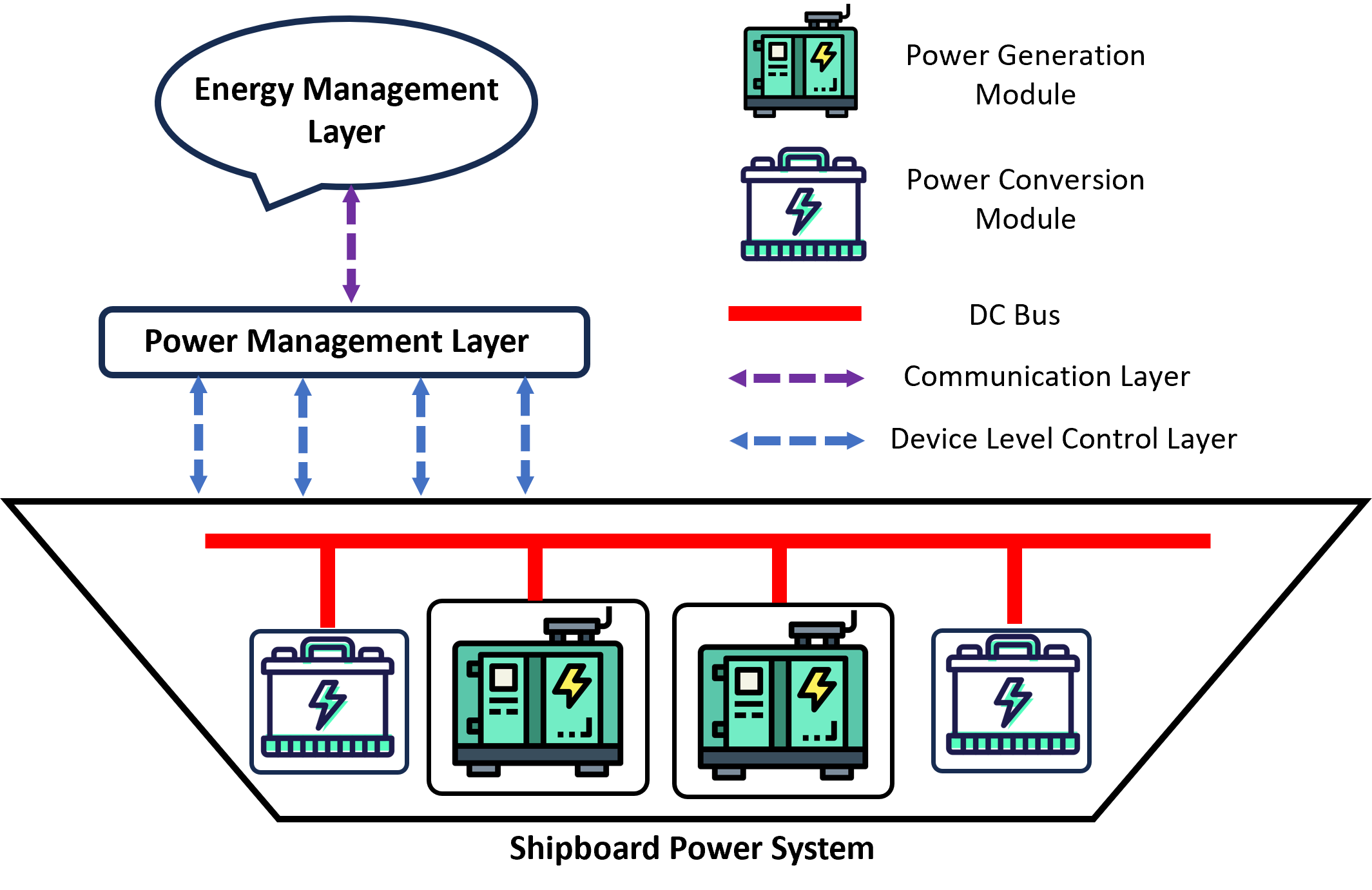}
	 \caption{An Islanded Notional Shipboard Power System}
     \label{4_Zone_SPS} 
\end{figure}

Existing literature considers the PGM and PCM ramp rate limitations in designing the energy management (EM) strategy \cite{2017_Vu_3,2021_Vedula,2017_Zohrabi,Zhang_2022}, which makes use of the optimal control theory (OCT) techniques such as dynamic programming (DP) and model predictive control (MPC). While limitations in the computational capabilities are the major roadblock in using the MPC over DP, it offers a better environment to incorporate constraints such as ramp rate limitations, and lower and upper power limitations on the PGMs and PCMs. Thus, due to recent advancements in the computational capabilities and its competence in handling the constraints, MPC, which originated as a control strategy for process control, made its way into the energy management of power systems \cite{Rawlings_book}\cite{PS_MPC_book}.

Nonetheless, integrating PCMs into the existing SPS framework raises another important objective of managing the degradation of the PCMs. It is a well-known fact that the PCMs degrade faster than the traditional gas-turbine-based PGMs. Thus, the EM problem needs to address the \textit{PCM degradation, PGM efficient operation and the ramp rate limitations} together. Since PCM degradation is a complex process that is difficult to capture mathematically with high precision, investigators have come up with a few degradation minimization measure heuristics embedded into the cost function. Minimizing these heuristics impacts the degradation process. Table-\ref{tab:references} shows the various degradation measures used as heuristics in the existing literature. 

\begin{table}[ht]
\centering
\caption{PCM Degradation Measure Heuristics and references}
\label{tab:references}
\resizebox{0.8\columnwidth}{!}{\begin{tabular}{c|c}
\hline \hline
\textbf{Degradation Measure }    & \textbf{References}  \\
 \hline \hline
Depth of Discharge (DoD)    & \cite{Hein_2021}\cite{Steen_2021}\cite{Li_2023}\cite{Zhao_2023}   \\
State of Charge (SoC)   & \cite{Nawaz_2023}\cite{Wang_2022}\cite{Vedula_2} \\
Charge/Discharge Cycle  & \cite{Ji_2023} \\
\hline
\end{tabular}}
\end{table}

While, \cite{Hein_2021,Steen_2021,Li_2023,Zhao_2023,Nawaz_2023,Ji_2023,Wang_2022} have used different measures for capturing the degradation process, we consider the \emph{PCM power, PCM state of charge and PGM power} together in an optimization-based degradation measure heuristic framework to analyze and study the PCM degradation process. A centralized model-predictive energy management problem is formulated. The impact of the PCM degradation measurement heuristics chosen with objectives: minimizing the \textit{absolute power extracted} out of the PCM and minimizing the \textit{state of charge} of the PCM. The results compare which of the two heuristics better aids in mitigating the PCM capacity loss.

The rest of the paper is organized as follows: Section \ref{sec:notations} presents the mathematical notations and the preliminaries on shipboard power system. In Section \ref{sec: model}, mathematical formulation and the device level control (DLC) design are given. In Section \ref{sec: control development}, the model predictive energy management problem formulation is provided. Finally, a numerical simulation is presented as case studies using a \emph{single PGM, single PCM, and PLM}. The results and trade-offs are also discussed in Section \ref{sec: simulation}. Conclusions are presented in Section \ref{sec: conclusion}.

\section{NOTATIONS AND PRELIMINARIES}\label{sec:notations}

The set of natural numbers, real numbers, and positive real numbers is represented by $\mathbb{N}$, $\mathbb{R}$ and $\mathbb{R}_+$. A matrix with $n$ rows and $m$ columns is denoted as $\mathbb{R}^{n \times m}$. Lower-case alphabets denote scalars (for example $x \in \mathbb{R}$ and $x \in \mathbb{N}$). Lower-case bold alphabets denote the vectors (for example $\mathbf{x} \in \mathbb{R}^{n}$ and $\mathbf{x} \in \mathbb{N}^n$). $\underline{\mathbf{0}}$ and $\mathbf{1}$ denote the vector of zeros and ones.
 For any vector $\mathbf{x} \in \mathbb{R}^n$, $\|\mathbf{x}\|_2 \triangleq \sqrt{\mathbf{x}^\top\mathbf{x}}$ and $\|\mathbf{x}\|_1\triangleq\sum_{i=1}^{n}|\mathbf{x}_i|$, represents the 2-norm and the 1-norm, respectively (where $|.|$ denotes absolute value). The symbol $\preceq$ denotes the component-wise inequality i.e. $\mathbf{x} \preceq \mathbf{y}$ is equivalent to $\mathbf{x}_i \leq \mathbf{y}_i$ for $i=1,2,\hdots,n$.

 The concept of MPC is based on the optimal control theory and the \textit{Bellman's principle of optimality}. \textit{``Given the optimal sequence $u^*(0),\hdots,u^*(N-1)$ and the corresponding optimal trajectory $x^*(0),\hdots,x^*(N-1)$, the sub-sequence $u^*(k),\hdots,u^*(N-1)$ is optimal for the sub-problem in the horizon $[k,N-1]$ starting from optimal state $x^*(k)$ \cite{borrelli_bemporad_morari_2017}"}. 

 \subsection{Preliminaries: Notional Shipboard Power System}
 \begin{figure}[t!] 
	\centering
	\includegraphics[width=0.48\textwidth]{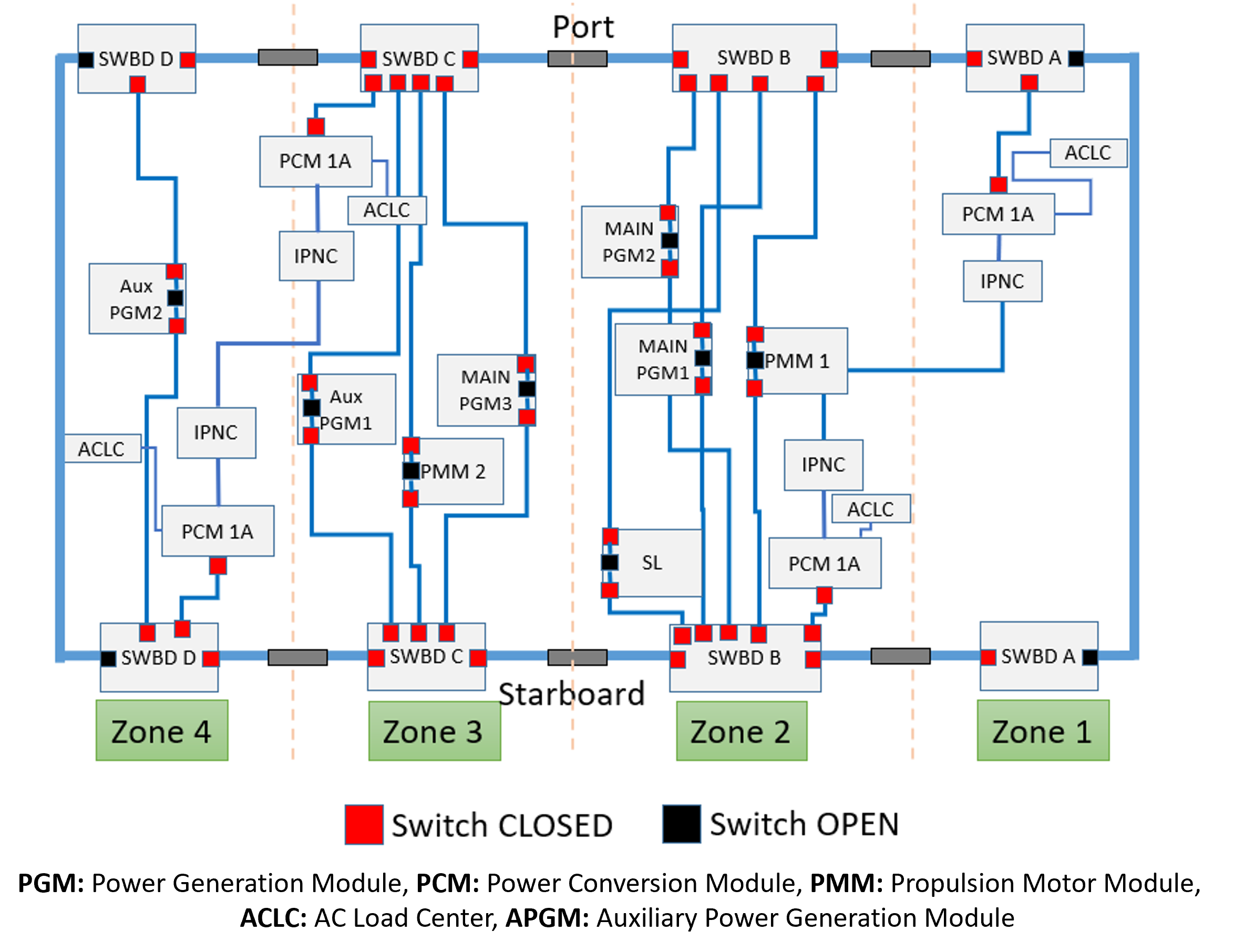} 
	\caption{Notional 4-Zone Shipboard Power System Model Presented by the Office of Naval Research \cite{ESRDC_1}. 
    }
	\label{SPS_Model}
\end{figure}
The shipboard power system model proposed by the U.S. Office of Naval Research (ONR) shown in Figure \ref{SPS_Model} comprises PGMs, PCMs, and device level controllers (DLCs) for local voltage and current control. These components are categorized into a zonal structure \cite{ESRDC_1}. PGMs consist of fuel-operated gen-sets. PCMs consist of multiple energy storage elements such as ultracapacitors and battery energy storage systems. All zones are unified through a common 12\textsf{kV} DC bus. The power demand must be met under any given circumstance i.e. power supplied by all sources must match the load demand scenarios such as high ramp and uncertainty in the load to maintain the system stability. In this paper, only the main power-supplying sources as mentioned before are considered. 
There are two main PGMs and two main PCMs. The combined power supplied by these sources must meet the SPS power demand. The following notations are used in the Shipboard Power System context:
\begin{itemize}
    \item PGM: Power Generation Module. It is also referred to as a generator in some parts of the paper. The usage of PGM and generator are synonymous.
    \item PCM: Power Conversion Module. It is also referred to as an energy storage system (ESS) or battery in some parts of the paper. The usage of PCM/ESS/Battery is synonymous.
    \item The assumption is that the ESS operates bidirectionally representing charging and discharging operations. Also, we assume that the ESS can ramp up faster than the generators.
    \item PLM: Power Load Module.
\end{itemize}

\section{Model Development}\label{sec: model}
We consider the lumped PGM and the lumped PCM models for the model development. Thus, the total number of the components in an SPS is reduced to a single PGM, a single PCM supplying a common single PLM. 

\subsection{Power Generation Module Modeling}
The PGM model used in this work is modeled as a controllable voltage source in series with a \textsf{RL}-line and a shunt capacitance. The mathematical model for the individual PGM is given as follows \cite{2021_Vedula}:
\begin{align}
    l_g \frac{di_g}{dt} &= -r_g i_g(t) + v_g - v_c, \\
    c_g \frac{dv_c}{dt} &= i_g - i_{g_{ref}},
\end{align}
where $i_g \in \mathbb{R}$ is the PGM current, $l_g, r_g, c_g\in \mathbb{R}_+$ are the inductance, resistance, and the capacitance of the PGM. The capacitance voltage/bus voltage $v_{c} \in \mathbb{R}$ and the controlled generator voltage $v_g \in \mathbb{R}$. Given a power reference ${p}_{g_{ref}} \in \mathbb{R}$ from the high-level control, the current reference for the device level PGM current controller (DLC) is given as ${i}_{g_{ref}} = {p}_{g_{ref}}/{v}_c \in \mathbb{R}$. The goal of the DLC is to achieve the objective:
$$\int_0^{\infty} \bigg(\underbrace{v_g(t)-v_{g_{ref}}(t)}_{\tilde{v}_g}\bigg)^2 dt < \infty.$$
Since the design and analysis of such a DLC is not the main objective of this paper numerous linear and nonlinear methods can be employed to design it \cite{Khalil_Book}. Thus, a simple \emph{PI} controller is employed of the form:
$$v_g = v_c + K_P \tilde{v}_g(t)+K_I \int_0^{t}\tilde{v}_g(\tau)d\tau,$$
where $K_P,K_I \in \mathbb{R}$ are the device level current control gains.

\subsection{Power Conversion Module Modeling}
The PCM is modeled as an ESS which comprises a single or multiple hybrid energy storage systems. This might include a collection of flywheels, supercapacitors, or battery energy storage systems (BESS). The dynamics of grid-following PCM used in this work and the battery current calculations are based on \cite{2021_Vedula}. The relationship between the PCM power injected $p_{b}$ and the SoC for the individual PCM is given as follows:
\begin{equation}\label{SoC_Power}
    \dot{q} = -\frac{p_b}{Q_bv_c}
\end{equation}
where $q \in \mathbb{R}$ denotes the State of Charge (SoC), ${Q}_{b} \in \mathbb{R}_+$ is the total capacity of the PCM in \textsf{AHr}. The bus voltage to which the PCM is coupled is given by $v_{c} \in \mathbb{R}$. The initial SoC is denoted by $q_0 \in [0,1]$. For this work, the SoC dynamics in (\ref{SoC_Power}) is discretized with a sampling time of ${T_s} \in \mathbb{R}_+$ using the Euler method, and the following discrete dynamics are obtained:

\begin{equation}\label{SoC_Power_Disctrete}
    q_{t+1} = q_t - \frac{T_s}{Q_b v_c}p_{b_t}.
\end{equation}

The PCM degradation model used in this work is an Arrhenius equation-based model which uses $Ah$-throughput as $\displaystyle \int_{0}^{t}\left|{i}_b(\nu)\right|d\nu$ as a metric to evaluate battery state of health (SoH). ${i}_b(t) \in \mathbb{R}$ is the current drawn from the PCM (positive while discharging and negative while charging). The PCM \textit{capacity loss} formulation is given as follows \cite{SONG2018433}.
\begin{equation}
    Q_L(t) =  e^{\frac{-\zeta_1+T_bC_{r}}{RT_b}}\int_{0}^{t}\left|{i}_b(\nu)\right|d\nu
\end{equation}
where $T_b$ is the PCM temperature in \textsf{Kelvin}, $C_r$ is the C-rate of the PCM. The capacity loss (\%) is given as follows: $$\Delta Q \% = \frac{Q_b-Q_{L}(t)}{Q_b}\times 100,$$
where $Q_{L}$ is the \emph{capacity loss} of the PCM in \textsf{ampere-hour}. \textit{This does not capture the battery end of life (EOL)}. It only captures the capacity lost during the PCM operation ($\Delta Q$).

\subsection{Power Load Module Modeling}
The PLM is a controllable load and is modeled as a current sink regulated through a controllable voltage. The dynamics are given as follows:
\begin{equation}
    l_L \frac{di_L}{dt} = -r_L i_L + \underbrace{v_L - v_c}_{\tilde{v}}.
\end{equation}

The goal of the PLM-DLC is to achieve the objective:
$$\int_0^{\infty} \bigg(\underbrace{i_L(t)-i_{L_{ref}}(t)}_{\tilde{i}_L}\bigg)^2 dt < \infty.$$
Given a load power $p_L$, the load current reference can be determined as $i_{L_{ref}} = p_L/v_c$. A simple \emph{PI} controller is employed of the form:
$$\tilde{v} = K_P \tilde{i}_g(t)+K_I \int_0^{t}\tilde{i}_g(\tau)d\tau,$$
where $K_{P_L},K_{I_L }\in \mathbb{R}$ are the PLM device level current control gains. 

Thus, the power flow in the SPS is given as follows:
\begin{equation}\label{SPS_Power_Flow}
    {{p}}_{g}(t)+{p}_{{b}}(t) - {p}_{{L}}(t) = 0,
\end{equation}
where ${p}_{g_i}, {p}_{b_i}, {p}_L\in\mathbb{R}$ are powers corresponding to the PGMs, the PCMs, and the PLMs respectively.

\begin{figure}[t!] 
	\centering
	\includegraphics[width=0.42\textwidth]{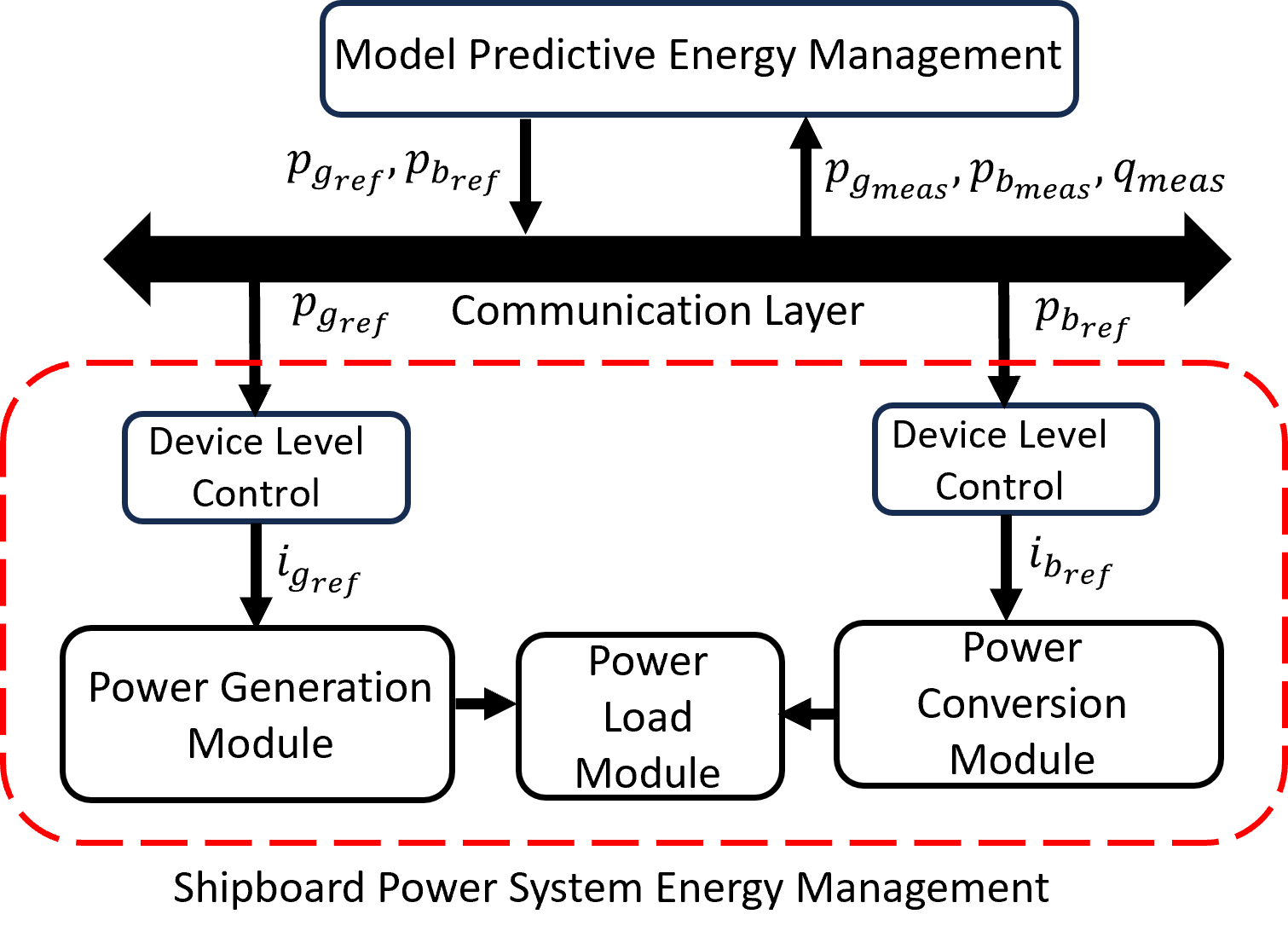} 
	\caption{Hierarchical control structure employed for the Energy Management in SPS}
	\label{EM_SPS_pic}
\end{figure}

\section{MODEL PREDICTIVE ENERGY MANAGEMENT}\label{sec: control development}
In this section, the optimization constraints imposed are developed followed by the model predictive energy management optimization problems considering \textit{battery power} and \textit{battery state of charge} as the PCM degradation heuristic measures.
\subsection{Constraints Development}
\textit{Ramp Limitations:}
Since the PGM is ramp-limited and the PCM has faster ramping capabilities compared to the PGM, the following constraints are imposed limiting the power-delivering capabilities in consequent time instances in a horizon ($H$).
\begin{subequations}\label{ramp-limitations}
    \begin{align}
        \left|\mathbf{p}_{g_k} - \mathbf{p}_{g_{k-1}}\right| \preceq r_g, \hspace{2mm} \forall k=1,\hdots,H, \\
        \left|\mathbf{p}_{b_k} - \mathbf{p}_{b_{k-1}}\right| \preceq r_b, \hspace{2mm} \forall k=1,\hdots,H,
    \end{align}
\end{subequations}
where $r_b > r_g$ and $r_g, r_b \in \mathbb{R}$ are the ramp-rate limitations of the PGM and the PCM.

\textit{Safety Constraints:} In order to accommodate for the safe operation of the PGM and the PCM, the lower and upper power limitations, SoC are imposed as
\begin{subequations}\label{safety-cons}
    \begin{align}
        \underline{\mathbf{p}}_{g} \preceq \mathbf{p}_{g_k} \preceq \overline{\mathbf{p}}_g, \hspace{2mm} \forall k=1,\hdots,H, \\
        \underline{\mathbf{p}}_{b} \preceq \mathbf{p}_{b_k} \preceq \overline{\mathbf{p}}_b, \hspace{2mm} \forall k=1,\hdots,H, \\
        \underline{\mathbf{q}} \preceq \mathbf{q}_k \preceq \overline{\mathbf{q}}, \hspace{2mm} \forall k=1,\hdots,H.
    \end{align}
\end{subequations}

\textit{SoC Constraint:} The discretized SoC dynamics in (\ref{SoC_Power_Disctrete}) are imposed as an equality constraint to be satisfied at every instant of the horizon:
\begin{equation}\label{soc-balance}
    \mathbf{q}_{k+1} = \mathbf{q}_k - \frac{T_s}{Q_b v_c}\mathbf{p}_{b_k}, \hspace{2mm} \forall k=1,\hdots,H.
\end{equation}

\textit{Power Balance Constraint:} The power supply and the demand must be satisfied in a grid at any given time instant. To that extent, the following constraint is imposed based on (\ref{SPS_Power_Flow}). It is assumed that the load power measured is held constant for the length of the optimization horizon interval.
\begin{equation}\label{power-balance}
    \mathbf{p}_{g_k} + \mathbf{p}_{b_k} - p_L\mathbf{1} = 0, \hspace{2mm} \forall k=1,\hdots,H.
\end{equation}

\subsection{Energy Management Problem Formulation}
The objective/cost of the energy management problem consists of three optimization variables: PGM power, PCM power, and SoC. The cost is formulated to minimize the PGM power injection, PCM power injection, and PCM SoC. Since the device level control is designed and is tasked to track the power set points provided by the MPC layer, the MPC problem and the imposed constraints are given as
\begin{equation}\label{MPEM_prb}
    \begin{aligned}
        \Minimize_{\mathbf{p}_g,\mathbf{p}_b,\mathbf{q}} \quad & \frac{\beta}{2}\norm{\mathbf{p}_{g}-\mathbf{p}_{g_r}}_2^2 + \frac{\gamma_p}{2}\norm{\mathbf{p}_{b}}_2^2 + \frac{\gamma_q}{2}\norm{\mathbf{q}-\mathbf{q}_0}_2^2 \\
        \SubjectTo \quad & (\ref{ramp-limitations})-(\ref{power-balance}),
    \end{aligned}
\end{equation}
where $\mathbf{p}_{g_r} \in \mathbb{R}^H$ is the known desired operating point of the PGM. $\mathbf{q}_0 \in \mathbb{R}^H$ is the initial SoC of the PCM. The objective is to maintain the PGM around a desired operating point and to maintain the SoC as close as possible to the initial SoC. The constraint set is a \textit{polytope} comprised of an affine equality constraint and an affine inequality constraints. The objective function is strongly convex. Based on this fact the optimization problem in (\ref{MPEM_prb}) is feasible, the optimization problem is convex whose global optimal solution can be found in polynomial time by existing algorithms \cite{boyd_vandenberghe_2004}.

\section{Numerical Simulation}\label{sec: simulation}
\begin{figure}[t!] 
	\centering
	\includegraphics[width=0.42\textwidth]{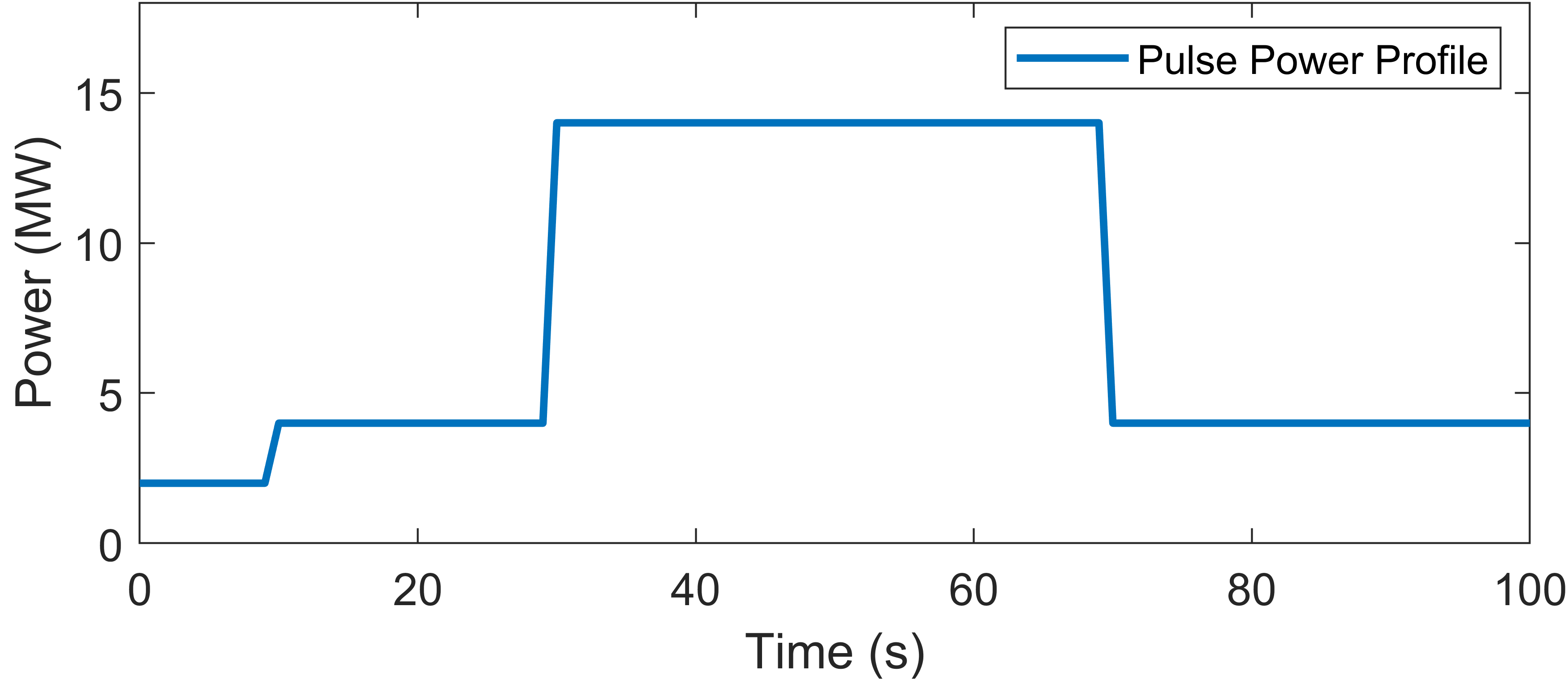} 
	\caption{Power Profile}
	\label{Power_Profile}
\end{figure}
A lumped model of the SPS consisting of a single PGM, a single PCM, and a single PLM is considered for simulation purposes. The PGM and PCM combined to supply the PLM. The developed model and the optimization problems are implemented on the Digital Storm desktop with the processor Intel(R) Core(TM) i9-14900K 3.20GHz and 64GB of RAM.  The simulation time-step was chosen to be 1\textsf{ms}. The time taken by the algorithm to converge was observed to be 0.04\textsf{seconds} and about 50 iterations for the optimal values to be within the tolerance limits. Thus, the MPC and model rate transition in Simulink simulation was set at 1\textsf{second} i.e. the MPC updates every second. The MPC was implemented using YALMIP \cite{lofberg}. The power profile is designed based on the PLM ratings and an occurrence of a pulse is considered from 20-70sec. The simulation parameters are shown in Table \ref{tab:rated}. Figure \ref{Power_Profile} shows the designed PLM load profile. Three simulation scenarios are considered based on the optimization problem (\ref{MPEM_prb}).
\begin{itemize}
    \item \textit{No PCM Heuristic (Scenario 1):} $\beta = 1$, $\gamma_p=0$, and $\gamma_q = 0$.
    \item \textit{Power Minimization Heuristic (Scenario 2):} $\beta = 1$, $\gamma_p=1000$, and $\gamma_q = 0$.
    \item \textit{SoC Minimization Heuristic (Scenario 3):} $\beta = 1$, $\gamma_p=0$, and $\gamma_q = 1000$.
\end{itemize}
The choice of the penalty parameters $\beta$, $\gamma_p$, and $\gamma_q$ is chosen based on the observations from the simulations undertaken with different penalty values. The chosen values for the penalty best capture the results of the optimization problem for different scenarios. 

\begin{figure}[h!] 
	\centering
	\includegraphics[width=0.40\textwidth]{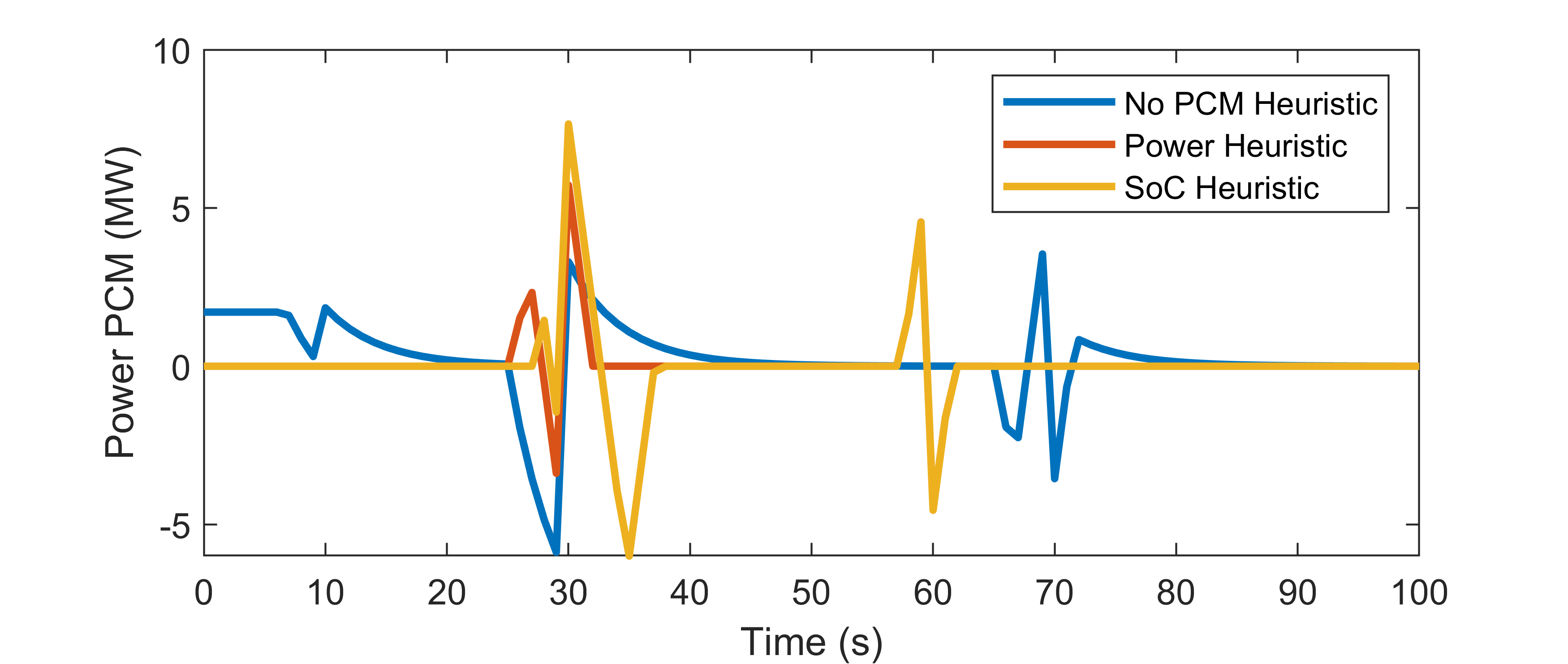} 
	\caption{PCM Powers for different Heuristic Scenarios}
	\label{Power_PCM}
\end{figure}

Figure \ref{Power_PCM} shows the PCM powers with different scenarios implemented. It can be observed that Scenarios 2 and 3 display almost similar trends in terms of power injected by the PCM in response to the requested load profile. Scenario 2 and 3 which penalizes the power extracted and the SoC deviation from the PCM provides the most conservative response to the requested load. 

\begin{table}[ht]
\centering
\caption{Rated Values and Simulation Parameters \cite{ESRDC_1}}
\label{tab:rated}
\resizebox{0.8\columnwidth}{!}{\begin{tabular}{c|c|c}
\hline \hline
\textbf{Parameter}  & \textbf{Parameter}  & \textbf{Parameter}  \\
\textbf{Description} & \textbf{Notation}   & \textbf{Value}  \\  \hline \hline
PGM Upper Power Limit   & $\overline{p}_g$   & 28 MW \\ 
PGM Lower Power Limit   & $\underline{p}_g$   & 0.2 MW \\ 
PCM Upper Power Limit & $\overline{p}_b$ & 10 MW \\
PCM Lower Power Limit & $\underline{p}_b$ & -10 MW \\
SoC Lower Limit & $\underline{q}_b$ & 0.7 \\ 
SoC Upper Limit & $\overline{q}_b$ & 0.8 \\
Horizon & $H$ & 5 sec \\ 
PCM Capacity (Li-ion) & $Q_b$   & 20 AHr \\ 
PCM Ramp-rate (Li-ion) & $r_b$   & $\overline{p}_b$ \\
PGM Ramp-rate  & $r_g$   & 0.1$\overline{p}_g$ \\\hline
\end{tabular}}
\end{table}

\begin{figure}[h!] 
	\centering
	\includegraphics[width=0.42\textwidth]{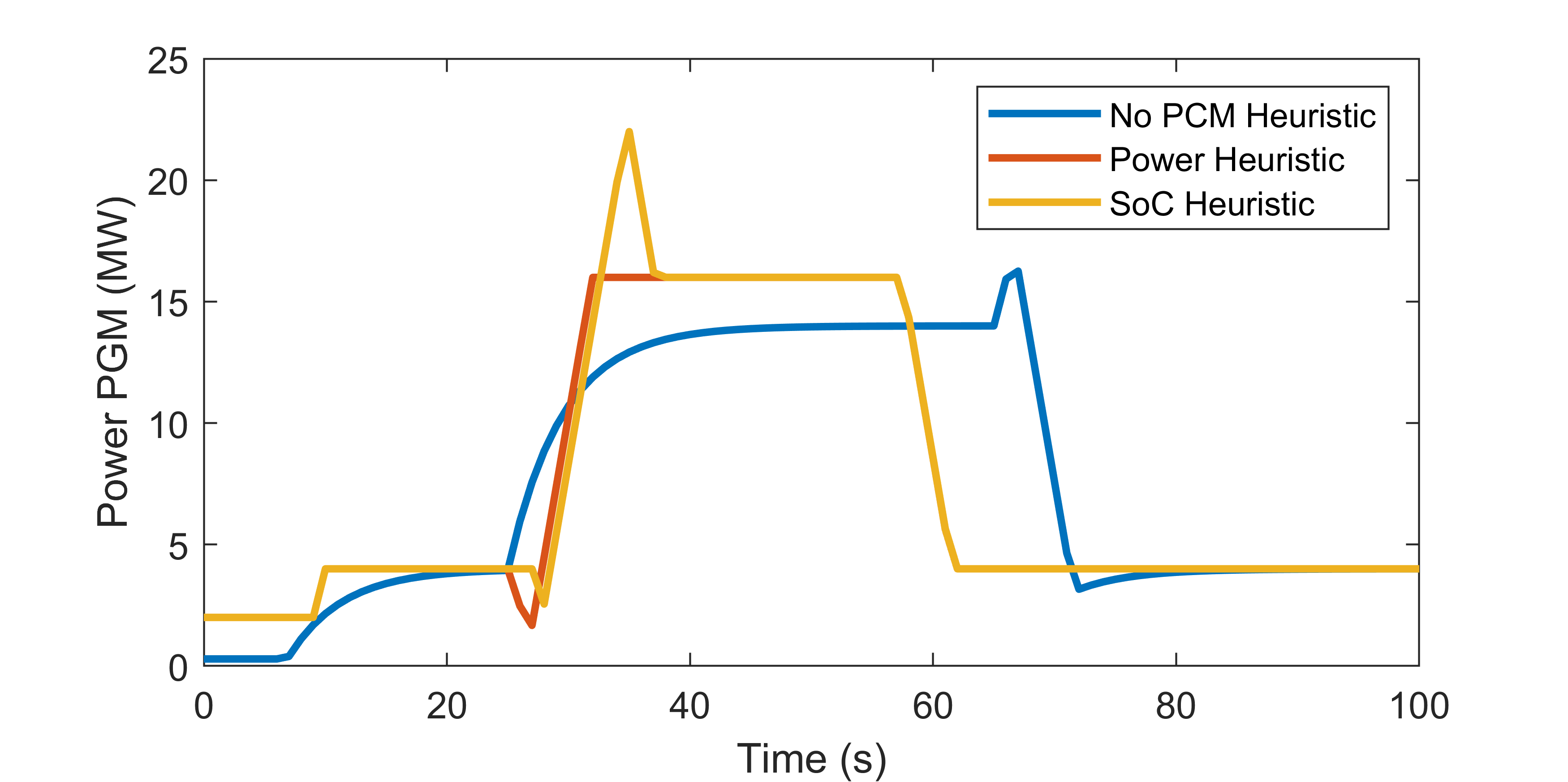} 
	\caption{PGM Powers for different Heuristic Scenarios}
	\label{Power_PGM}
\end{figure}

The impact of the optimization scenarios and the PCM minimization heuristics on the PGM can be seen in Figure \ref{Power_PGM}. The PGM is pushed to its ramping and operational limits when the PCM power and the SoC are minimized i.e. Scenario2 and Scenario 3 as compared to Scenarios 1. The PGM acts in response to maintain the power balance constraint.

\begin{figure}[h!] 
	\centering
	\includegraphics[width=0.42\textwidth]{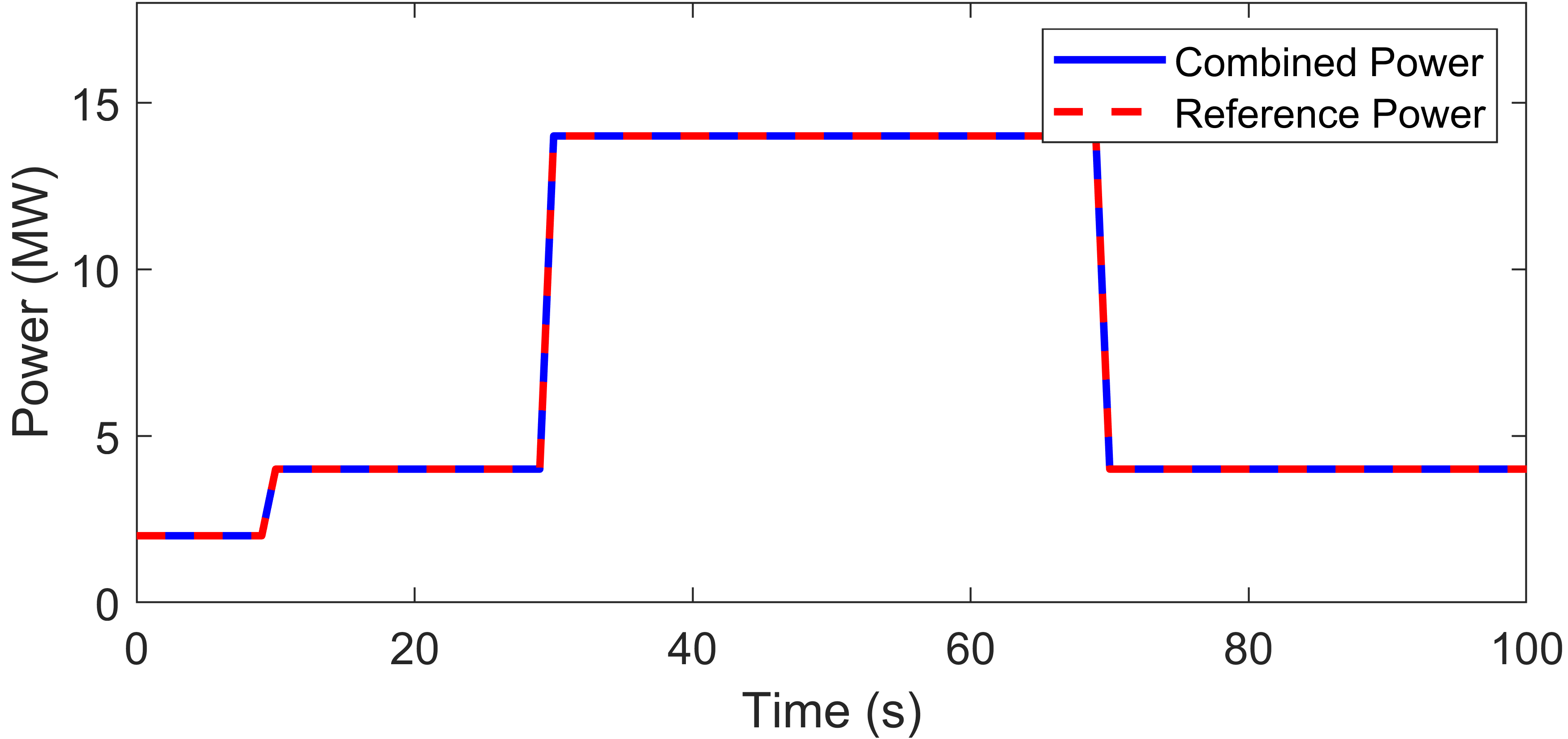} 
	\caption{Power Tracking response for Scenario 2}
	\label{Power_Tracking}
\end{figure}

The combined powers injected by the PGM and the PCM in response to the designed pulse power profile are shown in Figure \ref{Power_Tracking}. In order to avoid redundancy, the tracking response is plotted for Scenario 2. A similar response in tracking was also observed for other optimization scenarios.

\begin{figure}[h!] 
	\centering
	\includegraphics[width=0.45\textwidth]{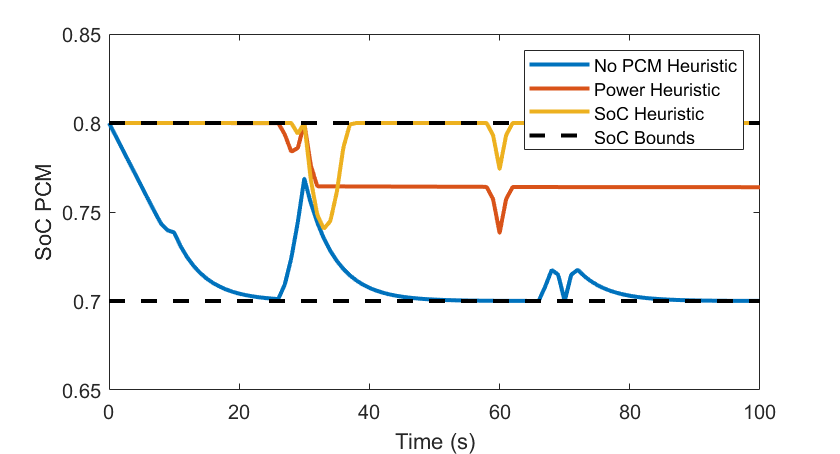} 
	\caption{SoC of the PCM for different Heuristic Scenarios}
	\label{PCM_SOC}
\end{figure}

Figure \ref{PCM_SOC} shows the State of Charge response for the proposed optimization scenarios. It is observed that all the scenarios obey the imposed upper and lower constraints on the SoC. Scenario-3 whose objective is to minimize the SoC, tries to keep the SoC as close as possible to the initial SoC. The SoC that is closest to the initial SoC can be observed during Scenario 3. 

\begin{figure}[h!] 
	\centering
	\includegraphics[width=0.42\textwidth]{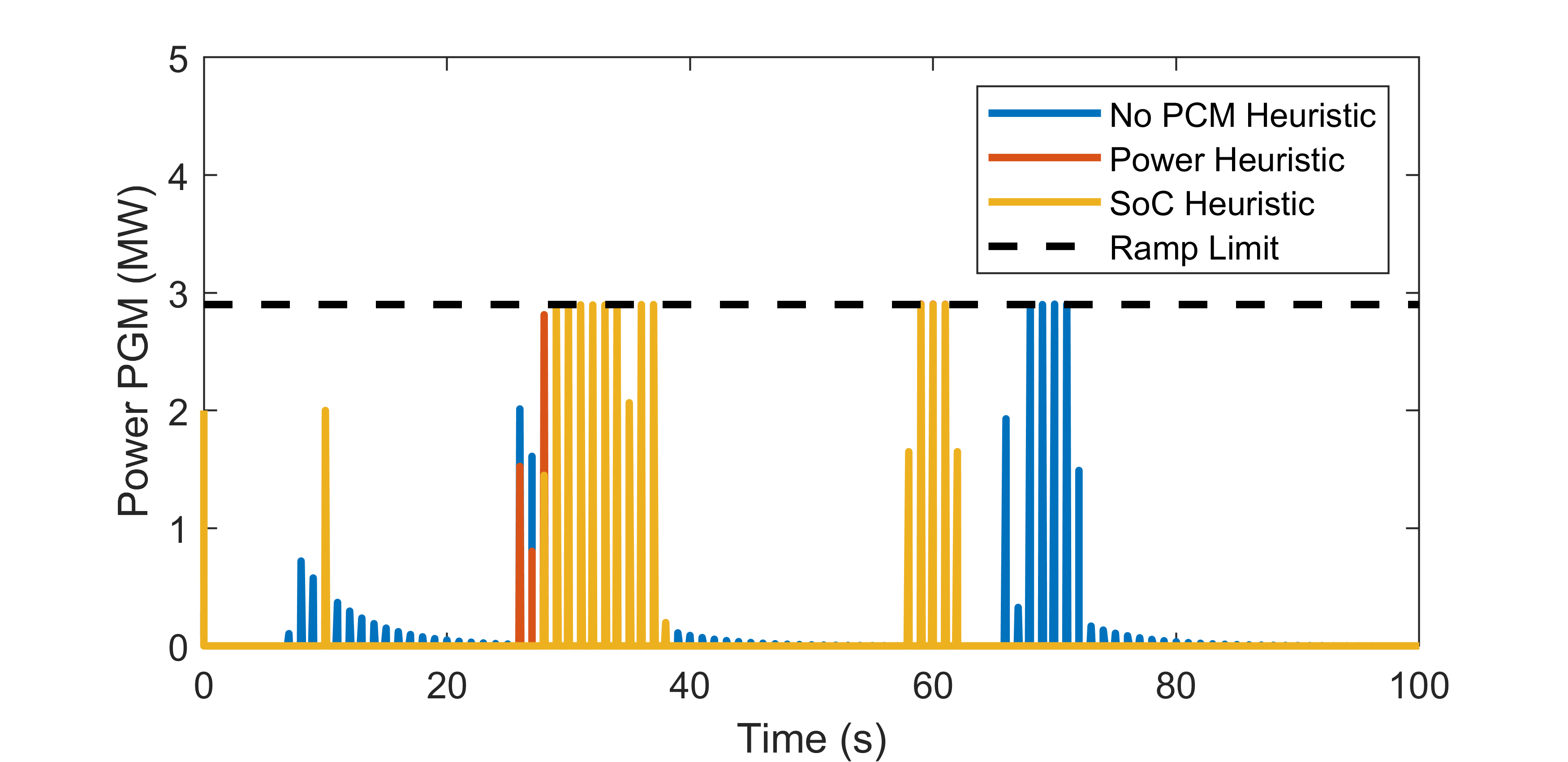} 
	\caption{PGM Ramp limits: PGM Power for every Horizon step}
	\label{PGM_ramp}
\end{figure}

\begin{figure}[t!] 
	\centering
	\includegraphics[width=0.42\textwidth]{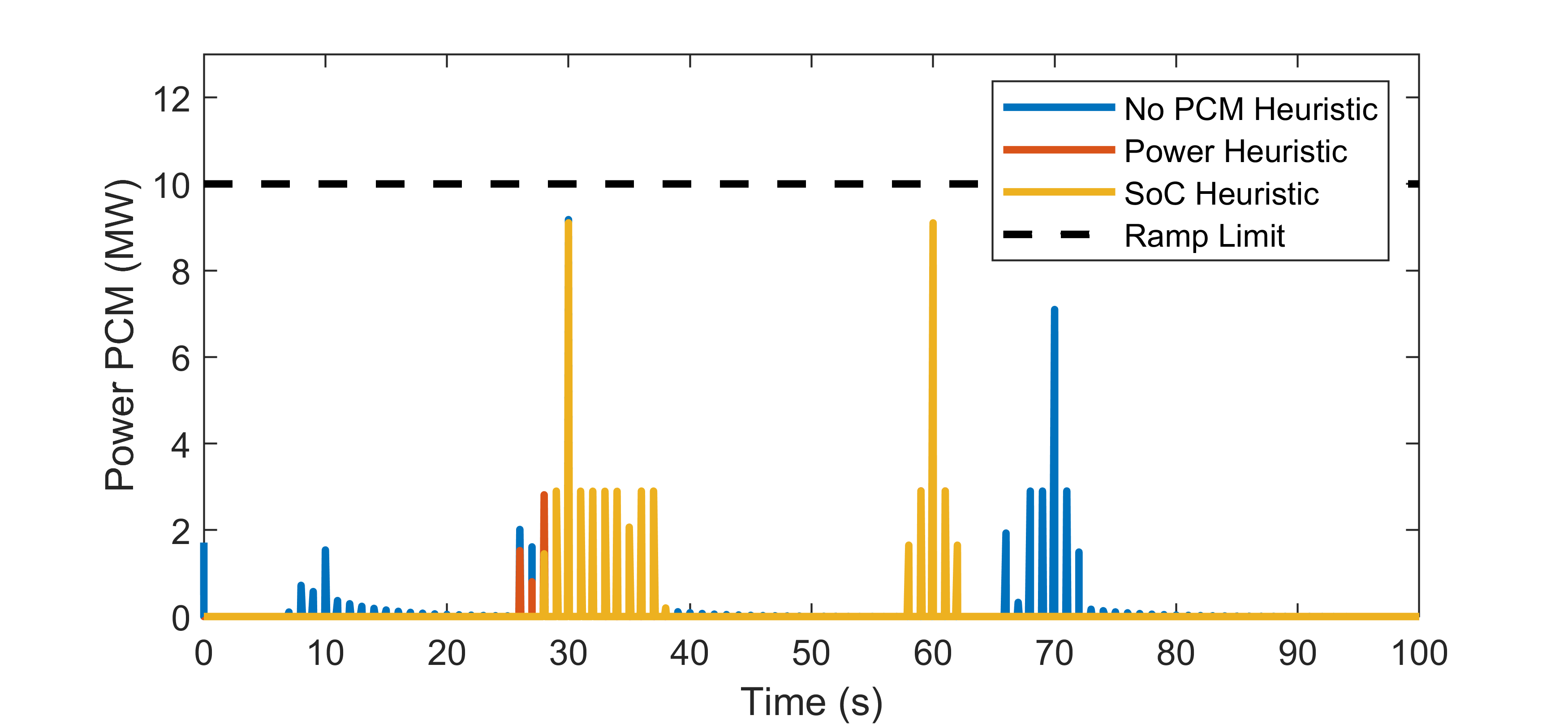} 
	\caption{PCM Ramp limits: PCM Power for every Horizon step}
	\label{PCM_ramp}
\end{figure}

The ramp limits for the PGM and the PCM are shown in Figure \ref{PGM_ramp} and Figure \ref{PCM_ramp}. It can be seen that Scenario 2 and Scenario 3 pushes PGM to the maximum ramp limit, while for the PCM the same scenario has the lowest ramping. This indicates a more conservative PCM usage at the expense of aggressive PGM usage in order to satisfy the power balance constraint. Another key observation during the ramping of PGM is around 70sec, where both Scenarios 2 and 3 exhibit similar PGM power injections. 

\begin{figure}[t!] 
	\centering
	\includegraphics[width=0.42\textwidth]{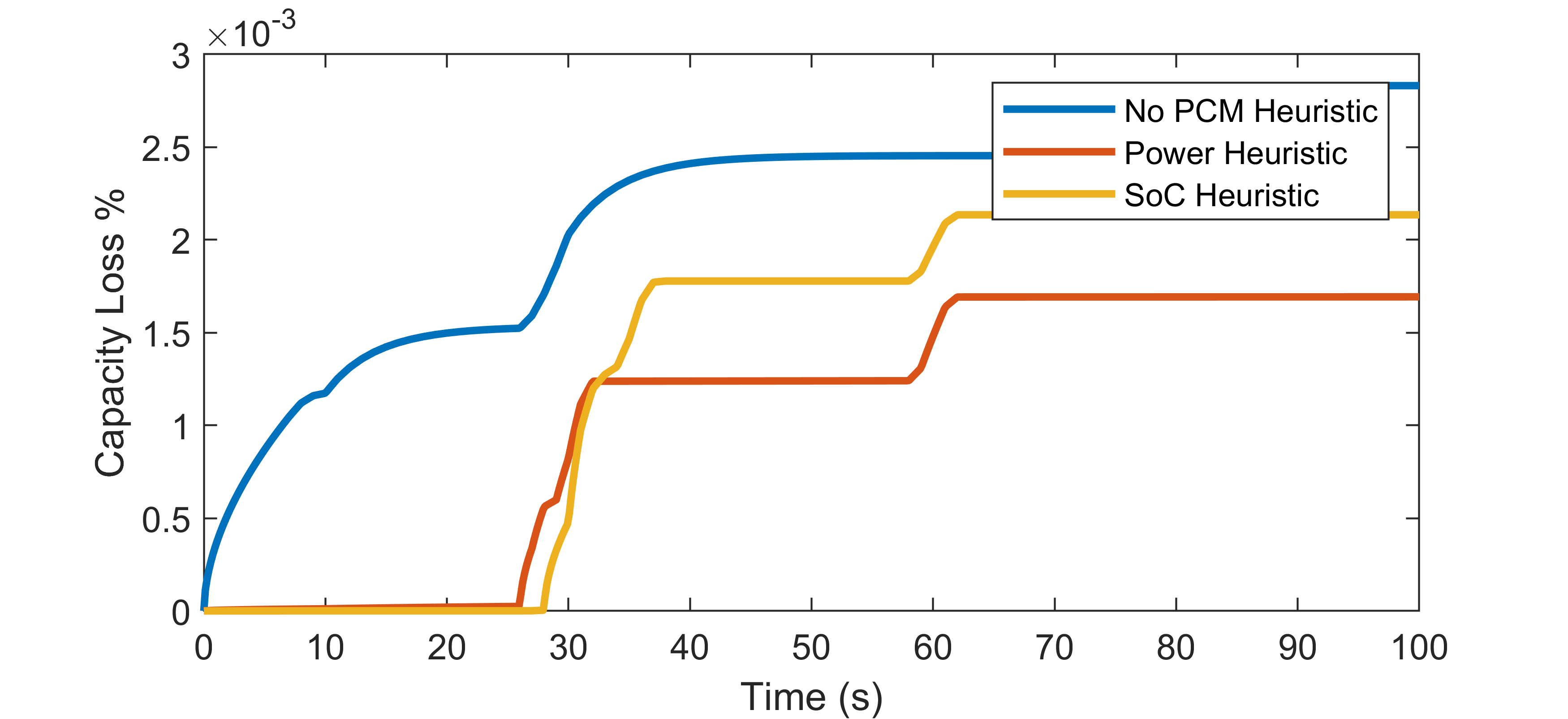} 
	\caption{Capacity Loss $\%$ for different Heuristic Scenarios}
	\label{Cap_Loss}
\end{figure}

Finally, Figure \ref{Cap_Loss} shows the \textit{Capacity Loss $\%$} of the PCM during its operation for different scenarios. It can be seen that Scenario 2 has a minimal impact on the PCM capacity loss for the operation duration, especially during the occurrence of pulse power. 

Thus, from the results, and testing different Scenarios it can be deduced that the PCM degradation minimization comes at the expense of pushing the PGM to its operational limits. However, this trade-off can be tuned by adjusting the penalties $\beta,\gamma_p$, and $\gamma_q$.

\section{Conclusion}\label{sec: conclusion}
This paper presents an MPC-based battery (PCM)-degradation-aware energy management strategy for a shipboard power system under high-ramp-rate pulse power load situations. The individual power generation module, power conversion module ramp rate limitations, the PGM-rated condition, and multiple PCM degradation heuristics in terms of PCM power and PCM SoC are considered in the problem formulation while allocating the power in the SPS. The model-based PCM degradation measures based on the PCM power and PCM state of charge are used to capture the PCM usage and thus minimize it. The effect of the PCM power minimization and PCM state of charge has been observed on the PCM capacity loss. Finally, a numerical example has been given to show the effectiveness of the proposed EM method and study which heuristics have a greater impact in mitigating PCM degradation.

\end{document}